\documentclass[11pt]{article}
\usepackage[latin1]{inputenc}
\usepackage[T1]{fontenc}
\usepackage[english]{babel}
\usepackage{amssymb}
\usepackage{amsmath} % pour quelques options mathematiques
\usepackage{amsfonts} % pour utiliser les polices ams
\usepackage{geometry} % pour changer la taille des marges
\usepackage{graphicx} % idem
\usepackage[hyperindex,colorlinks=false]{hyperref} % pour ajouter des hyperliens
\usepackage{color} % pour ajouter de la couleur
\usepackage{subfigure,comment}
\newcommand{ \V }[1]{ \underline{#1} }
\newcommand{ \M }[1]{ \underline{\underline{#1}} }
\newcommand{\gradot}{{}^{t}\underline{\underline{\nabla}}_\textit0}
\newcommand{\grado}{\underline{\underline{\nabla}}_\textit0}

\newcommand{\fd}{{\underline{f}_d}}

\newcommand{\Omegao}{{\Omega_\textit0}}

\newcommand{\GammaEo}{{\Gamma_{E_\textit0}}}
\newcommand{\GammaEEpo}{{\Gamma_{E^{\phantom{'}}_\textit0E'_\textit0}}}

\newcommand{\Gammao}{{\Gamma_\textit0}}
\newcommand{\E}{{E^{\phantom{'}}_\textit0}}

\newcommand{\OmegaEo}{{\Omega_{E_\textit0}}}

\newcommand{\dOmegaEo}{{\partial \Omega_{E_\textit0}}}
\newcommand{\uE}{{\underline{u}_{E_\textit0}}}

\newcommand{\uEstar}{{\underline{u}^{\star}_{E_\textit0}}}

\newcommand{\EpsilonpEstar}{{\underline{\underline{\dot{E}}}(\underline{u}^{\star}_{E_\textit0})}}
\newcommand{\GLE}{{\underline{\underline{E}}_{E_\textit0}}}

\newcommand{\piE}{{\underline{\underline{\pi}}_{E_\textit0}}}

\newcommand{\WE}{{\underline{W}_{E^{\phantom'}_\textit0}}}
\newcommand{\WEM}{{\underline{W}_{E_\textit0}^M}}
\newcommand{\WEMM}{{\underline{W}^M_{E_0}}}
\newcommand{\WEm}{{\underline{W}_{E_\textit0}^m}}
\newcommand{\WEstar}{{{\underline{W}^{\star}_{E_\textit0}}}}

\newcommand{\FEo}{{\underline{F}_{E^{\phantom'}_\textit0}}}

\newcommand{\FEoM}{{\underline{F}^M_{E_\textit0}}}
\newcommand{\FEom}{{\underline{F}^m_{E_\textit0}}}
\newcommand{\FEpo}{{\underline{F}_{E_\textit0'}}}

\newcommand{\ktEo}{\left({k^{m}_{E_\textit0}}\right)_t}
\newcommand{\knEo}{\left({k^{m}_{E_\textit0}}\right)_n}
\newcommand{\kmEo}{{k^-_{E_\textit0}}}

\newcommand{\kpEo}{{k^+_{E_\textit0}}}

\newcommand{\sWE}{{\mathcal{W}_{E_0}}}
\newcommand{\sWEO}{{\mathcal{W}_{E_0}^0}}

\newcommand{\sWME}{{\mathcal{W}_{E_0}^M}}

\newcommand{\sFE}{{\mathcal{F}_{E_0}}}

\newcommand{\sFME}{{\mathcal{F}_{E_0}^M}}
\newcommand{\sAd}{\mathbf{A_d}}
\newcommand{\sGamma}{{\boldsymbol{\Gamma}}}

\newcommand{\suEO}{{\mathcal{U}_{E_0}^0}}

\newcommand{\sEp}{{\mathbf{E^+}}}
\newcommand{\sEm}{{\mathbf{E^-}}}

\newcommand{\FchapEo}{{\underline{\widehat{F}}_ {E_\textit0}}}

\newcommand{\WchapE}{{\underline{\widehat{W}}_{E_\textit0}}}

\title{Virtual delamination testing through non-linear multi-scale computational methods: some recent progress}

\author{O. Allix$^{1}$, P. Gosselet$^{1}$, P. Kerfriden$^{2}$, K. Saavedra $^{3}$\\ 
\textit{\small{$^1$ LMT-Cachan, ENS Cachan/CNRS}} \\
\textit{\small{61, avenue du pr\'esident Wilson, F-94230 Cachan, France}}\\
\textit{ \small{$^2$ now at Cardiff University, School of Engineering,  }} \\
\textit{ \small{Queen's Buildings, The Parade, Cardiff CF24 3AA, UK}}  \\
\textit{ \small{$^3$ now at Departamento de Tecnolog\'ias Industriales, Universidad de Talca, Los Niches km 1, Curic\'o, Chile }}
}

\begin{document}
\maketitle

\abstract{This paper deals with the parallel simulation of delamination problems at the meso-scale by means of multi-scale methods, the aim being the Virtual Delamination Testing of Composite parts. In the non-linear context, Domain Decomposition Methods are mainly used as a solver for the tangent problem to be solved at each iteration of a Newton-Raphson algorithm.  In case of strongly non linear and heterogeneous problems, this procedure may lead to severe difficulties. The paper focuses on methods to circumvent these problems, which can now be expressed using a relatively general framework, even though the different ingredients of the strategy have emerged separately. We rely here on the micro-macro framework proposed in \cite{Ladeveze01b}. The method proposed in this paper introduces three additional features: (i) the adaptation of the macro-basis to situations where classical homogenization does not provide a good preconditioner, (ii) the use of non-linear relocalization to decrease the number of global problems to be solved in the case of unevenly distributed non-linearities, (iii) the adaptation of the approximation of the local Schur complement which governs the convergence of the proposed iterative technique. 
Computations of delamination and delamination-buckling interaction with contact on potentially large delaminated areas are used to illustrate those aspects.}

\section{Introduction}
%Despite the many studies since the beginning of the eighties, the prediction of delamination remains a challenge from a scientific as well as an industrial point of view.  Many difficulties are inherent to composite modeling and characterization and the interested reader can find a survey on the large bibliography devoted to the modeling of composites in \cite{Herakovich12}. In this framework the simulation of delamination leads do deal with specific difficulties  (i) the complex state of stress which leads to the initiation and the propagation of delamination (ii) the size and the complexity of the problem to be solved if one wish to deal with composite structures. One of the industrial objectives is to replace some of the numerous tests performed today in order to assess the damage tolerance of composites by numerical simulations \cite{Allix06}, in other words to perform Virtual Delamination Testing (VDT).

Despite the many studies since the beginning of the eighties, the prediction of delamination remains a challenge, both from a scientific and from an industrial point of view. Many difficulties encountered in this field are inherent to composite modeling. The interested reader will find a survey of the literature devoted to the modeling of composites in \cite{Herakovich12}. Simulating delamination requires to tackle specific difficulties (i) the complex state of stress which leads to the initiation and the propagation of delamination (ii) the size and the complexity of the problem to be solved when dealing with composite structures. One of the industrial objectives is to replace some of the numerous physical tests required to assess the damage tolerance of composites by numerical simulations \cite{Allix06}, or in other words to perform Virtual Delamination Testing (VDT).

%It is today well accepted at industrial level that a realistic VDT procedure shall rely on a meso-scale description of the laminates. This implies an underlying Finite Element description whose characteristic size is at most the thickness of individual plies, roughly the tenth of millimeter. For the treatment of industrial parts, the size of the resulting model is therefore huge, leading to consider the applicability of parallel computations to those large and complex non-linear cases. Indeed these problems involve possible multiple delaminations, buckling and contact, and therefore instabilities.% Moreover quite a lot of questions are connected to the residual strength after impact, which implies, in turns to deal with combined delamination and buckling and post-buckling analyses.

It is today well accepted at an industrial level that a realistic VDT procedure should rely on a meso-scale description of the laminates. Such a fine-scale description requires a Finite Element discretisation whose characteristic size is at most the thickness of individual plies, which is of the order of a tenth of a millimetre. For the treatment of industrial parts, the size of the resulting model is therefore huge, which leads to considering the applicability of parallel computations to those large and complex non-linear cases. Indeed, these problems potentially involve multiple delamination fronts, buckling or contact, and therefore instabilities.

The paper describes some attempts to deal with such problems within a multi-scale parallel framework when delamination is the only damage mode that is taken into account. %The open issue regarding how to take into account, within a multi-scale computational framework,  the coupling between inner-layer damage mode, as transverse cracking, and delamination is discussed in the final part of the paper. 
Delamination is here modeled by means of a cohesive zone model. Several such models have been proposed in the literature, for example in \cite{Allix92,Xu94,Schellekens94,Alfano01b,Camanho02}. In this paper we make use of the model proposed in  \cite{Allix96c} but similar results could be obtained with other models once calibrated. The main features of the model are briefly presented in Section~\ref{sec:model}. 

%The  paper then focuses on some aspects of the underlying numerical multi-scale framework that we follow to address such question.  
In the remainder of the paper, the basics of the multi-scale strategy are described and test cases are presented to demonstrate its numerical efficiency. The contribution of the paper is the definition of a framework that encompasses different aspects of the strategy that have emerged separately in \cite{Ladeveze07,Guidault07,cresta07,kerfriden09,Allix2010,Saavedra12b}. Indeed, our previous publications focussed on particular issues related to the present topic. It now appears that the set of tools that we have developed in the past can help formulate a general non-linear multi-scale approach that is adapted to the treatment of large laminated and damageable structures taking into account possible buckling and contact.
%The main aspects of the underlying numerical multi-scale framework are then described.  The multi-scale strategy relies on previous works wherein different aspects have emerged separately.  The strategy can now be settled in a relatively general formulation of a non-linear multi-scale approach that is adapted to the treatment of large laminated and damageable structures. 
Multi-scale refers here to a computational technique which involves different scales and which aims at finding the solution over the whole structure at the smallest scale of interest.
In this domain, one of the main issues to be addressed is the transfer of information from one scale to another.

A large number of multiscale methods rely on classical homogenisation of heterogeneous media, which was initiated by \cite{Hill63}. For linear periodic media the most efficient method was  initiated in \cite{Sanchez-Palencia80}. Further developments and related computational approaches can be found in  \cite{Oden96,Fish97,Schrefler99,Kouznetsova02,Feyel03}. An overview of these methods can be found in \cite{ZW08}. 
In computational homogenization approaches the resolution of the macro-problem leads to effective values of the unknowns; then, the micro solution is calculated locally based on the macro solution. The fundamental assumption, besides periodicity, is that the ratio of the characteristic length of the small scale to the characteristic length of the large scale must be small. All areas which do not satisfy the hypothesis of scale separation such as  boundary zones or vicinity of cracks, where the material cannot be homogenized, require special treatment. Therefore a large amount of recent studies aim at applying enhanced homogenization schemes to localized failure, see for instance \cite{massartpeerlings2007,belytschkoloehnert2007,coenenskouznetsova2012}. These extension, which are dedicated to overcome some of the difficulties for a given type of problem are not suited yet to deal with edge effects and large cracks and therefore not directly suitable to the treatment of delamination. 

The proposed method is based on a Domain Decomposition approach. These methods have been developed initially for large linear problems for which they can now be considered mature. The most well-known methods are FETI \cite{farhat94} and BDD \cite{mandel93} which where designed for perfect interface and extended, in the case of FETI, to contact (with or without friction) \cite{dostal07} or prescribed displacement gaps \cite{kruis07}.
FETI and BDD rely on a static condensation for each sub-domain. The resulting interface problem, which is large, is solved iteratively using a preconditioner which is associated to the balance of the rigid-body momentum of the sub-domains. In the case of heterogeneous media, the macro-space of rigid body is not necessarily the most appropriate one and a specific deflation can be performed \cite{Efendiev11}. It is therefore interesting to define the macro-space as the space which is associated with the homogenized response of the subdomains. The use of homogenization in the context of a Domain Decomposition Method has been proposed in the micro-macro approach of \cite{Ladeveze99c}.

Homogenization being a natural tool for composites, we base our approach on this micro-macro framework. The decomposition of the structure into subdomains is chosen to be compatible with the mesomodeling of the laminate. 
%The method being a mixed one, enables then to introduce any kind of interfacial constitutive law at the level of the interfaces between the subdomains \cite{Ladeveze07}. 
The micro-macro method relies on mixed interface conditions between subdomains, which allows to treat any kind of interfacial constitutive law directly at the interface level \cite{Ladeveze07}.
Three types of interfaces are considered here:  damageable interfaces and contact interfaces that are compatible with the mesomodeling, and perfect interfaces within the plies to split the problem into smaller pieces that are suited to parallelism.
This leads to a large number of substructures, and therefore to a large macro-problem. A third scale is then introduced to perform a low cost and precise approximation of the homogeneized macro-problem by using a BDD preconditioner \cite{kerfriden09}, as described in Section~\ref{sec:mainfeatures}. Note that in this paper only unilateral contact is considered but the framework is well suited to friction \cite{Champaney-2008}.

In the context of non-linear problems, Domain Decomposition Methods are mainly used as a solver for the tangent problem to be solved at each iteration of a Newton-Raphson based algorithm \cite{letallec94,germain07b,Klawonn07}, with if required an arc-length control. Several problems can be encountered in this case.  One of such problems is that the choice of the macro-space, and the implied homogenized behavior, is in general not adapted to a situation where a crack interacts with a subdomain. In other words, a preconditioner designed for a problem without crack can become ineffective in the presence of a crack. It has been shown in \cite{Guidault07,Guidault08}, that for cracks crossing interfaces of the domain decomposition, a modification of the macro-basis that adds a few degrees of freedom only makes it possible to obtain the same convergence rate for a cracked domain and for an un-cracked one. 
Unfortunately such an adaptation is not possible in the case of cracks propagating at the interface between subdomains. We have therefore followed another line, the one of non-linear relocalization \cite{cresta07,pebrel08}. This approach aims at solving the problem of the strong deterioration of the convergence rate due to localized non-linear effect \cite{cai-2002}. Following this line of thoughts, we have applied a non-linear relocalization approach on the set of subdomains that are connected to the front of delamination, which is described in Section~\ref{sec:reloc} \cite{Allix2010}. Note that our approach naturally avoids to conduct useless nonlinear computations far from the delamination front which is the goal of the prediction procedure developped in \cite{lloberas12} for Newton-Schur strategies.

Having in mind the problem of the tolerance of laminates to compression after impact  \cite{Demoura97} we have started to work on the question of the interaction between potentially large growing delaminated areas and buckling, taking into account the contact condition \cite{Saavedra12b}. In order to reconnect the micro part of the solution, the iterative technique introduces parameters which are associated to the influence of neighboring subdomains.  These so called ``search directions'' correspond to Robin parameters whose values were optimized for 3D structures. These parameters must be adapted to the aspect ratios of the slender structures by introducing well-chosen anisotropic coefficients. Moreover, in the case of contact between slender structures over large delaminated areas, the search directions should be optimized according to the interface's state (open or closed) in order to obtain a reasonable convergence rate of the iterative solver. These points are discussed in Section~\ref{sec:adaptbuck}. However, changing these parameters requires a refactorization of the homogenized operators, which is computationally expensive. A compromise between convergence speed and unit cost of iterations must be found so that the total CPU time is minimised.

%Finally, because of the stiffness loss resulting from buckling and delamination, the macrostiffness and search directions may become irrelevant and the macroscopic operators may need to be adjusted in order for the homogenized behavior to better represent the current state of the structure, as illustrated in Section~6.

\section{The interface model}\label{sec:model}
\begin{figure}[htb]
       \centering
       \includegraphics[width=0.8 \linewidth]{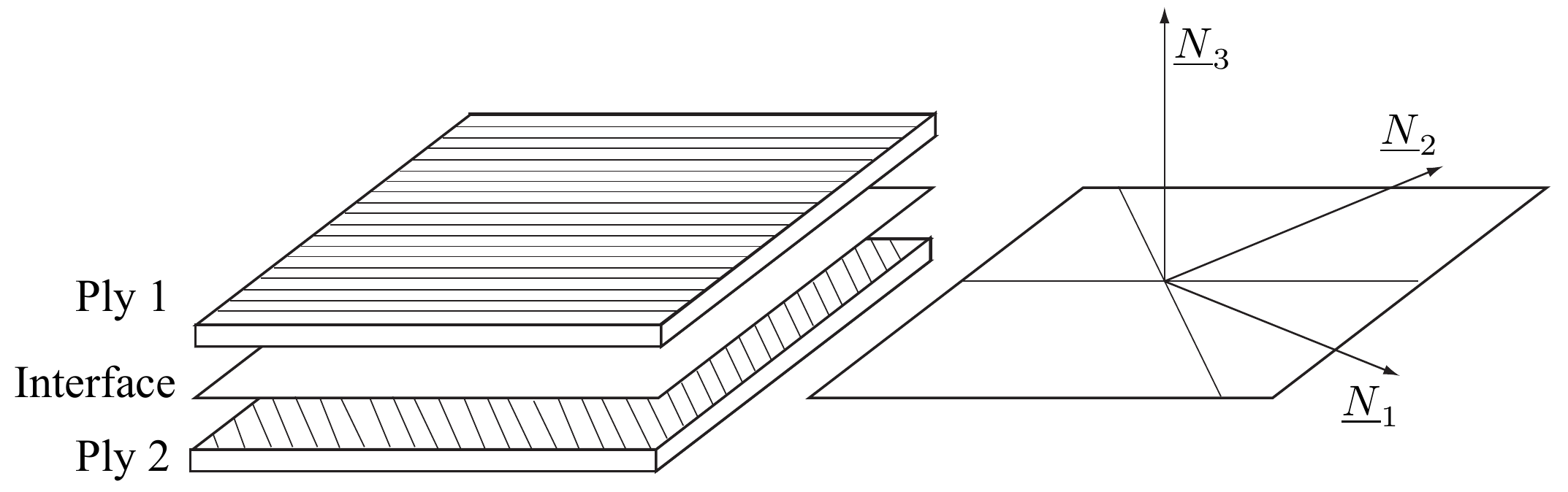}
       \caption{The mesomodel entities}
       \label{fig:interface_composite}
\end{figure}

The interface model relates the jump of displacements $\V{[u]} = \V{u}_{p'} - \V{u}_{p}$ to the normal Cauchy stress  $\V{t} = \M{\sigma}_{p}.\V{n} = \M{\sigma}_{p'}.\V{n}$ on the interface $\Gamma$ between two plies $p$ and $p'$. Here, $\V{n}$ denotes the outer normal to ply $p$ on $\Gamma_{pp'}$. In this problem, large displacements are considered but  the jump of displacement is assumed to be small prior to full delamination. The tractions $\V{t}$ are uniquely defined and related to the displacement discontinuities by means of the following damageable constitutive relation:
%\begin{equation}
%\left\{ \begin{array}{l}
%\displaystyle \M{\sigma}.\V{n} = \M{K}. \V{[u]} \\
%\displaystyle \M{\sigma}^+.\V{n} + \M{\sigma}^-.\V{n} = 0
%\end{array} \right.
%\end{equation}
\begin{equation}
\displaystyle \V{t} = \M{K}. \V{[u]} 
\end{equation}
The expression of the stiffness operator $\M{K}$ can be given in the reference frame $(\V{N}_1,\V{N}_2,\V{n})$ (as depicted in Figure (\ref{fig:interface_composite}) with $\V{n}=\V{N}_3$, $p \equiv \text{``Ply 1''}$ and $p' \equiv \text{``Ply 2''}$):
\begin{equation}
\left( \begin{array}{ccc}
\displaystyle (1-d_1){k_1} & 0 & 0 \\
0 & \displaystyle (1-d_2){k_2} & 0 \\
0 & 0 & \displaystyle \left(1- \left< \V{[u]}.\V{n} \right>_+ \ d_3 \right){k_3}
\end{array}
\right)
\end{equation}
$< \ >_+$ is here the positive indicator function. 
%The behavior equation $\M{\sigma}.\V{n} = \M{K}. \V{[u]}$ is equivalent to: $\M{\sigma}.\V{n} = \M{K}. \V{[u]}$

The local damage variables ${d_i}$ are introduced into the interface model in
order to simulate its softening behavior when the structure is
loaded. Their values range from $0$ (healthy
interface point) to $1$ (completely damaged interface point).
The parameters $d_i$ are related to the local energy release rates $Y_i$ of the interface's degradation modes as follows:

\begin{equation}
Y_i = - \frac{\partial e_d}{\partial d_i}
\qquad \textrm{where} \quad
\left\{ \begin{array}{ccl}
Y_1 & = & \displaystyle \frac{1}{2} k_1 [u]_1^2 \\
Y_2 & = & \displaystyle \frac{1}{2} k_2 [u]_2^2 \\
Y_3 & = & \displaystyle \frac{1}{2} k_{3} \left< [u]_3 \right>_+^2 \, ,
\end{array} \right. \end{equation}
where $e_d$ is the strain energy of the interface per unit area and $[u]_i$ is the $i^\text{th}$ component of field $\V{[u]}$. We assume that the damage variables are functions of a
single quantity: the maximum $Y_{|t}$ over time of a combination of
the energy release rates ${Y_i}_{|\tau}$, $\tau \leq t$:
\begin{equation}
Y_{|t} =  \underset{\tau \leq t}{\text{sup}} \left( {Y_3}_{|\tau}^\alpha + \gamma_1 {Y_1}_{|\tau}^\alpha + \gamma_2 {Y_2}_{|\tau}^\alpha \right)^{\frac{1}{\alpha}}
\end{equation}
%\begin{equation}
%\nonumber
%\textrm{where} \quad {Y}_{|\tau}  = \left( {Y_3}_{|\tau}^\alpha + \gamma_1 {Y_1}_{|\tau}^\alpha + \gamma_2 {Y_2}_{|\tau}^\alpha \right)^{\frac{1}{\alpha}}
%\end{equation}
Thus, the evolution law is defined by:
\begin{equation}
d_1 = d_2 = d_3 = w(Y) 
\end{equation}
\begin{equation}
\nonumber
\textrm{where, in general,} \quad w(Y) = \frac{n}{n+1} \left( \frac{Y}{Y_c} \right)^n\;.
\end{equation}

This model has the advantage, using a single damage variable, to be able to recover different critical energy release rates for pure mode loading as well as to fit to classical mixed mode criteria. For example, setting $\alpha=2$ leads to the classical quadratic energy criteria.

\section{Main features of the method}\label{sec:mainfeatures}

\begin{figure}[ht]
       \centering
       \includegraphics[width=0.99 \linewidth]{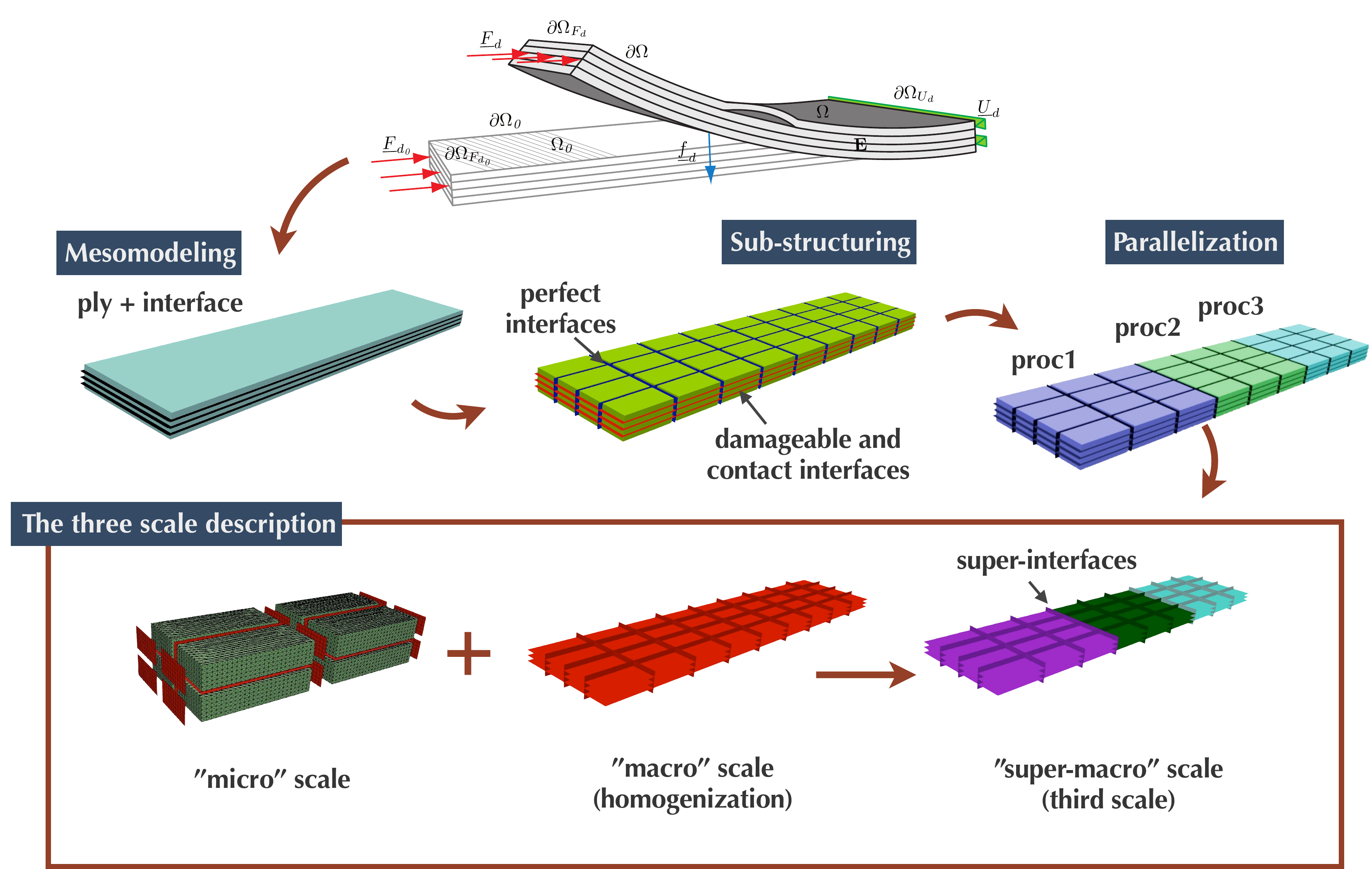}
       \caption{Decomposition of the problem and associated three computational scales}
       \label{fig:3scale}
\end{figure}

As presented in Figure~\ref{fig:3scale}, the structure is split into substructures which are connected by interfaces that have  mechanical behaviors. In a second stage, the substructures are gathered into ``super-substructures'' which are assigned to independent processors. In order to ensure the scalability of the method, the interfaces are the support of a macro (homogenized) problem, the resolution of which requires the solution of a third scale (supermacro) problem on super-interfaces.

\subsection{Micro-macro scale separation}
The iterative LaTIn algorithm, which is designed to non-linear problems, is here applied to the solution of the geometrically non-linear substructured reference problem with ``material'' non-linearities (damage, contact) localized on the interfaces. The finite element method is used to discretize the governing equations. The aim of the method is to find the subdomains fields $\uE$ (displacement), $\piE$ (Piola-Kirchhoff stress), and the interface fields $\WE$ (displacement), $\FEo$ (interforces), where index $E$ ranges over all substructures and index $0$ means that the quantity is given in the reference configuration (total Lagrangian approach).

In order to ensure the scalability of the method, a global and linear coarse grid problem is solved. The definition of the macroscopic fields required to construct this problem is done on the interface unknowns only. Whichever the interface behavior, action reaction principle $\FEo + \FEpo = \underline{0}$ is verified; the principle of the macro-problem is to enforce this equation partially at any time of the iterative solution process:
\begin{equation}\label{eq:macro_eq}
\forall {\WEMM}^\star\in \sWME,  \quad \int_{\GammaEEpo}  (\FEo + \FEpo)\cdot{\WEMM}^\star\, d \Gammao= 0 \, ,
\end{equation}
where displacement macrospace $\sWME$ is a parameter of the method as well as its dual $\sFME$.
The definition of the macro-spaces is done through a projector $\Pi$, which is chosen $L^2$-orthogonal (so that the same projector is used by the traction and displacement fields). We then have the following splitting of the interface quantities:

% These subspaces are common to neighboring substructures, they induce a splitting of interface quantities which is made unique through the uncoupling of virtual works:
% \begin{align}
% \displaystyle \FEo = \FEoM + \FEom \\ 
% \displaystyle \WE = \WEM + \WEm
% \end{align}
\begin{equation}
\begin{aligned}
&\forall  (\FEo,\WE)  \in \sFE \times\sWE, 
\\ &\FEo = \FEoM + \FEom \text{ with } \FEoM = \Pi \FEo  \text{ and }\WE = \WEM + \WEm \text{ with } \WEM = \Pi \WE 
\\
 &\Rightarrow\int\limits_{\GammaEEpo} \FEo\cdot\WE \, d \Gammao 
%\\ 
=  \int\limits_{\GammaEEpo} \FEoM\cdot\WEM \, d \Gammao  
  + \int\limits_{\GammaEEpo} \FEom\cdot\WEm \, d \Gammao
  \end{aligned}
\end{equation}

Numerical tests showed that in order to ensure the numerical scalability of the method the macro-basis should extract at least the linear part of the interface forces. Indeed, this macro-space contains the part of the interface fields with the highest wavelength. Consequently, according to the Saint-Venant principle, the micro complement only  has a local influence.

%%%%%%%%%%%%%%%%%%%%%%%%%%%%%%%
\subsection{Iterative algorithm}

It is now possible to split all the equations of the system  into two groups:
\begin{itemize}
    \item non-linear equations in the substructures and macroscopic admissibility of interfaces, whose solutions are elements of the manifold $\sAd$:
    \begin{enumerate}
        \item[-] non-linear kinematic admissibility of the substructures
        \begin{equation} \label{eq:kine_add_E1}
%\GLE=\frac{1}{2}\left(\FF_E^{t} \FFE - I_{d}\right) \;,\; \textrm{at each point of } \OmegaEo
\GLE= \frac{1}{2}\left(\grado \uE + \gradot \uE + \gradot \uE \gradot \uE\right) 
\ ,\ \text{on} \ \OmegaEo
\end{equation}
\begin{equation} \label{eq:kine_add_E2}
\uE_{| \dOmegaEo} = \WE_{|\GammaEo} \;,\; \text{on} \; \GammaEEpo
\end{equation}
        \item[-] non-linear static admissibility of the substructures
        \begin{multline} \label{eq:eq_sst}
\forall (\uEstar,\WEstar) \in \suEO \times\sWEO \quad \textit{such that} \quad \uEstar_{| \dOmegaEo} = \WEstar_{| \GammaEo},
\\
\int_\OmegaEo  \piE : \EpsilonpEstar \ d\Omegao 
\\
= \int_\OmegaEo \rho_{\E} \; \fd \cdot \uEstar \, d \Omegao + \int_{\GammaEo} \FEo \cdot \WEstar \, d \Gammao  
\end{multline}
        \item[-] behavior of the substructures
        \begin{equation} \label{eq:rdc_E}
\piE = \frac{\partial \psi}{\partial \GLE} \;,\; \text{on} \; \OmegaEo\;,\; %= \KEo \, \GLE (\uE)
\end{equation}
        \item[-] macroscopic admissibility of the interfaces \eqref{eq:macro_eq}.
    \end{enumerate}
    \item local (non-linear) equations in the interfaces whose solutions belong to the manifold $\sGamma$:
    \begin{itemize}
        \item[-] interface behavior (perfect, contact and cohesive interfaces) and boundary conditions (see Section~2).
    \end{itemize}
\end{itemize}

%The interface solutions $s = (s_\E)_{\E \in \struct} = (\WE , \FEo )_{\E \in \struct}$ of the first set of equations belong to the space $\sAd$, while the interface solutions $\widehat{s} = (\widehat{s}_\E)_{\E \in \struct} = (\WchapE , \FchapEo )_{\E \in \struct}$ of the second set of equations belong to $\sGamma$. 
The solution $s_{ref}=(\uE,\piE,\WE,\FEo)$ is such that:
\begin{equation}
s_{ref} \in \sAd \cap \sGamma
\end{equation}
Note that in the case of small perturbations, the admissibility equations are linear. In this case, the manifold $\sAd$ is an affine space.

\begin{figure}[ht]
       \centering
       \includegraphics[height=0.4 \textwidth]{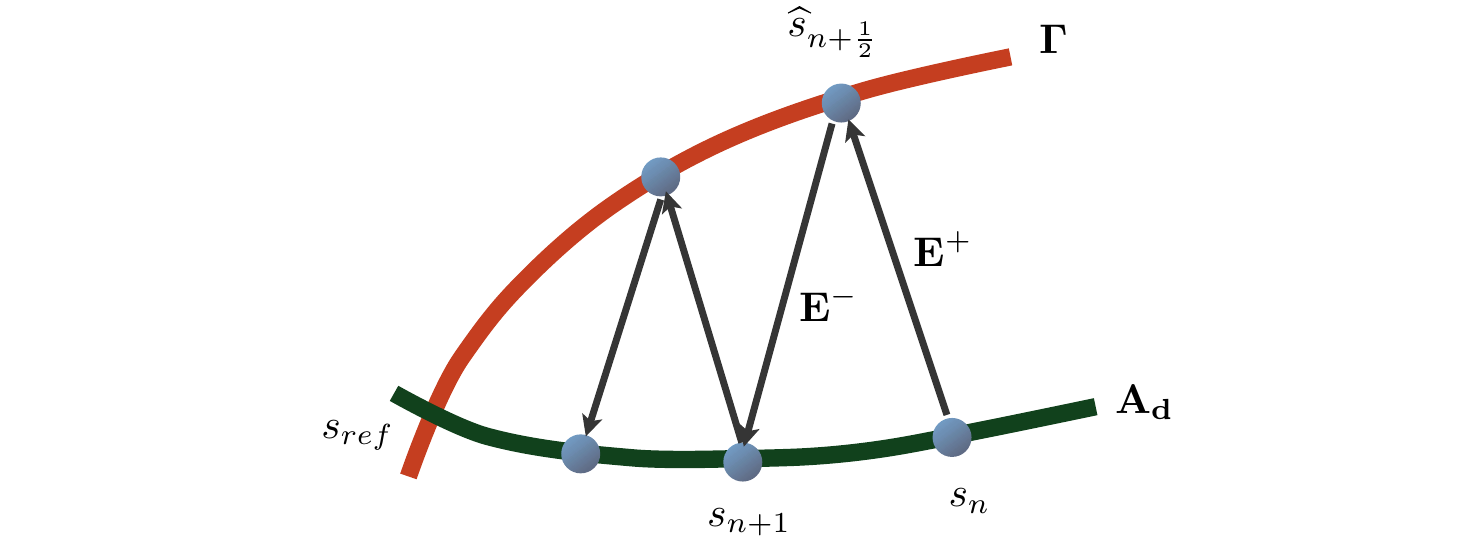}
       \caption{Schematic representation of the iterative LaTIn algorithm}
       \label{fig:algo_latin}
\end{figure}

The resolution scheme consists in seeking the solution $s_{ref}$ alternatively in these two manifolds: first, one finds a solution $s_n$ in $\sAd$, then a solution $\widehat{s}_{n+\frac{1}{2}}$ in $\sGamma$. In order for the two problems to be well-posed, search directions $\sEp$ and $\sEm$ linking the solutions $s$ and $\widehat{s}$ through the iterative process are introduced (see Figure \ref{fig:algo_latin}).
\begin{itemize}
\item At the global stage, starting from fields $(\WchapE,\FchapEo)$ satisfying the interface equations (manifold $\sGamma$), we seek fields $({\WE},{\FEo})$ in $\sAd$ (satisfying the subdomain equations) using the following closing relation:
\begin{equation}
\kmEo\left(\WE-\WchapE\right)+\left(\FEo-\FchapEo\right)=0
\end{equation}
$\kmEo$ is the search direction associated to the global stage.
\item At the local stage, starting from fields $({\WE},{\FEo})$ satisfying the subdomain equations (belonging to manifold $\sAd$), we seek fields $(\WchapE,\FchapEo)$ satisfying the interface equations (manifold $\sGamma$), using the following closing relation:
\begin{equation}
\kpEo\left(\WE-\WchapE\right)-\left(\FEo-\FchapEo\right)=0
\end{equation}
$\kpEo$ is the search direction associated to the local stage.
\end{itemize}

\subsection{Third scale}

\begin{figure}[htb]
       \centering
       \includegraphics[width=0.6 \linewidth]{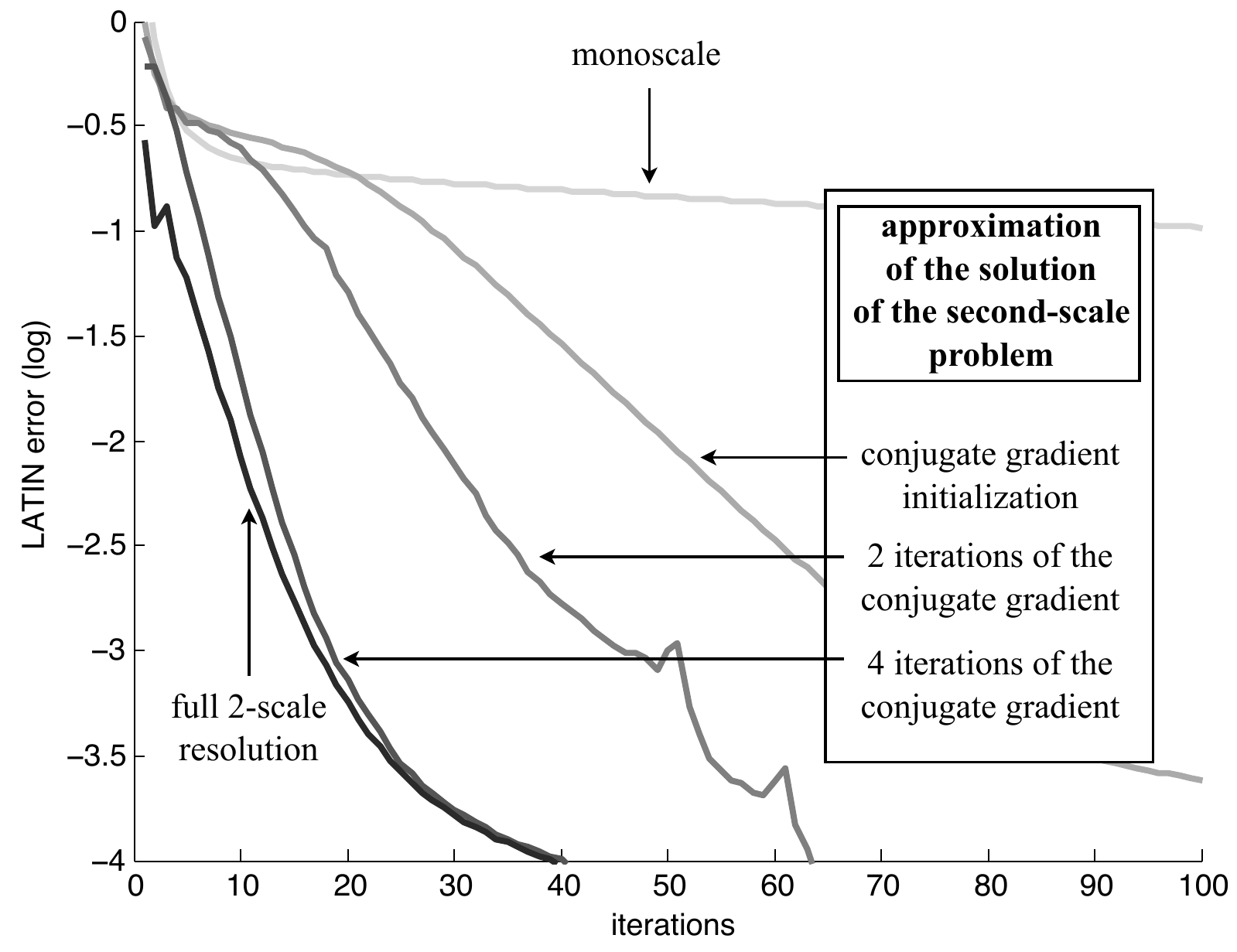}
       \caption{LaTIn convergence curves for different numbers of macro iterations}
       \label{fig:cg_conv}
\end{figure}
\begin{figure}[htb]
\centering
\includegraphics[width=0.8 \linewidth]{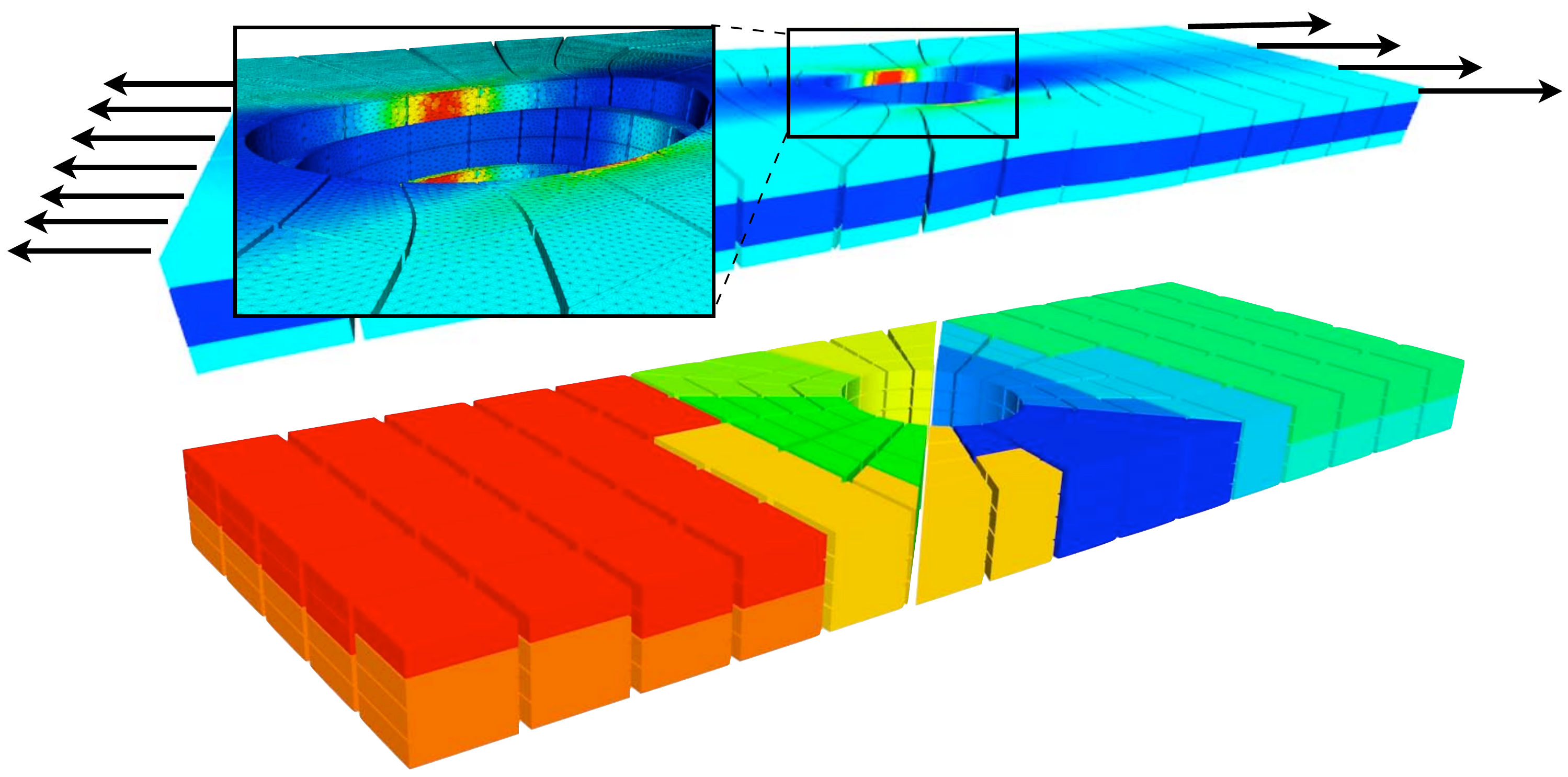}\caption{Different levels of substructuring}\label{fig:holed}
\end{figure}

The macro-problem being linear, discrete and sparse, it can be efficiently solved by a balancing domain decomposition algorithm \cite{mandel93}. To do so, substructures are gathered into super-substructures. The macro-problem is condensed at the super-interfaces and this condensed problem is solved by a conjugate gradient algorithm. The balance of super-structures requires the global solution of a third scale problem whose size is the number of rigid body motions of floating super-substructures.

Figure (\ref{fig:cg_conv}) shows the convergence rate of the LaTIn
algorithm when the conjugate gradient scheme for the condensed
macroproblem is stopped after a fixed number of iterations. The test
case is the holed plate under traction loading  with the super-substructuring pattern given
in Figure (\ref{fig:holed}). This problem is $3.4\ 10^6$ degrees of freedom large, it is distributed in 520 subdomains with $1,350$ interfaces leading to a $12\ 10^3$-unknown macro-problem. 8 super-substructures are employed so that the condensed macro-problem on the super-interfaces has $1,746$ unknowns and the dimension of the third scale (associated to the rigid body motions of the super-substructures) is only 36.

It appears clearly that only very few iterations of the conjugate gradient scheme are required to get the necessary part of the macro-displacement Lagrange multiplier leading to the multiscale convergence rate of the LaTIn algorithm. Typically, the algorithm is stopped when the residual error (normalized by the initial error) falls below $10^{-1}$.  The admissibility of the macro-forces is thus enforced on a third level, which is sufficient to determine the part of the solution which needs to be transmitted at each iteration of the resolution. In that case the time required to solve the macro-problem is divided by a factor 100.

\begin{figure}[ht]\centering
\includegraphics[width=.99\textwidth]{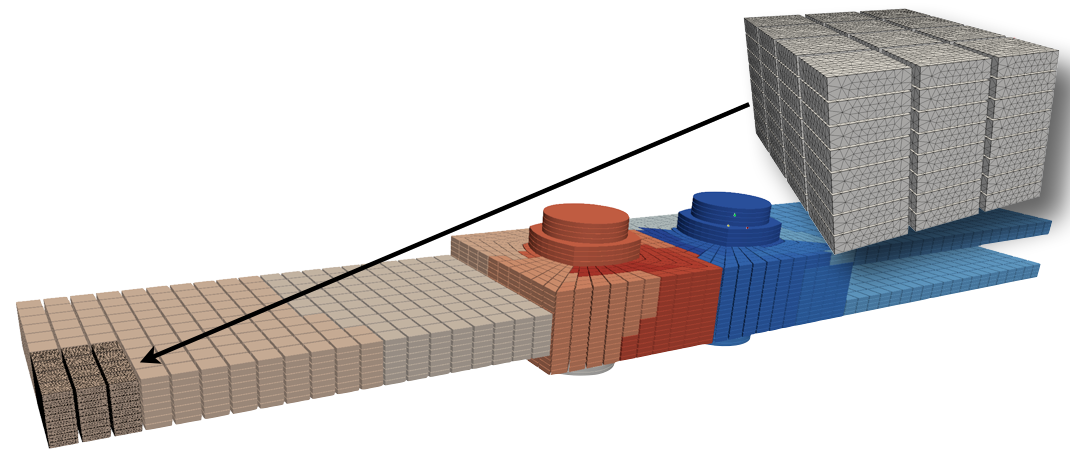}\caption{Bolted joint, substructuring and super-substructuring}\label{fig:liaisonsuper}
\end{figure}

\begin{figure}[ht]\centering
\includegraphics[width=.99\textwidth]{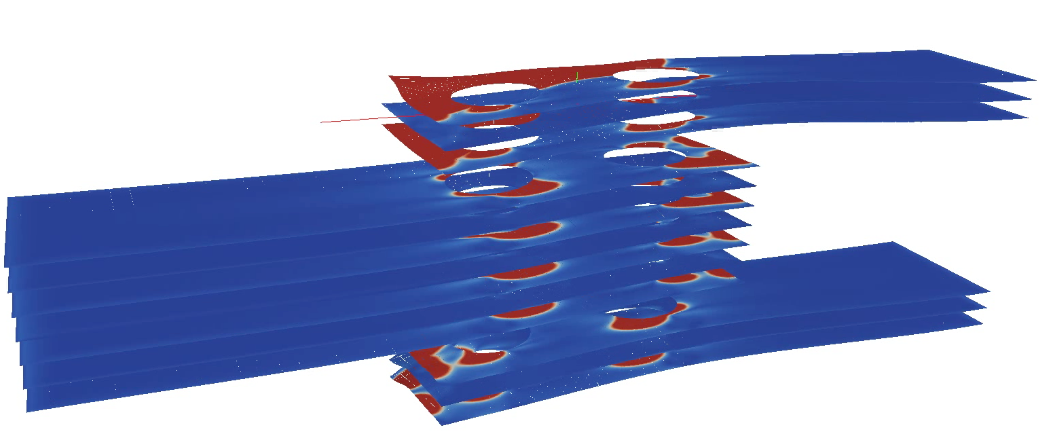}
\caption{Final delamination state of the bolted joint}\label{fig:delamliaison}
\end{figure}

Figure~\ref{fig:liaisonsuper} presents a larger case where a direct computation of the macroproblem would be too expensive. This  $12\ 10^6$ dof problem is distributed on $10,960$ substructures and $32\ 10^3$ interfaces, leading to a $300\ 10^3$ macro d.o.f, the super-substructuring leads to a $40\ 10^3$ condensed dof problems and $150$ super-macro d.o.f. Figure~\ref{fig:delamliaison} shows the extension of delamination in the final converged state.


%The previous examples allow us to make an important remark related to our vision of the problem of large-scale simulation of delamination. The reader will notice that we use a homogeneously fine finite element discretisation to approximate the displacement field in the structure. Specifically, the mesh is sufficiently fine to handle the propagation of a delamination front at an arbitrary location in the structure, which permits to handle multiple and arbitrary large delamination fronts without remeshing. The fact that we are not willing to remesh is a choice driven by the following observations (i) error estimation and adaptivity procedures are, in our opinion, not robust enough to allow for VDT, although some pioneer investigations can be found in \cite{ladevezemoes1999a,huertarodriguezferran2002,larssonrunesson2003,allixkerfriden2010}, and (ii) remeshing implies the development of field transfer strategies, which is relatively difficult in the context of softening nonlinear behaviors. Of course, the price that we pay for the expected robustness associated to this choice is a large number of arguably unnecessary degrees of freedom.

The previous examples allow us to make an important remark related to our vision of the problem of large-scale simulation of delamination. The reader will notice that we use a homogeneously fine finite element discretization to approximate the displacement field in the structure. Specifically, the mesh is sufficiently fine to handle the propagation of a delamination front at an arbitrary location in the structure, which permits to handle multiple and arbitrary large delamination fronts without remeshing. The fact that we are not willing to remesh is a choice driven by the observation that error estimation, mesh adaptivity and field transfer are, for parallel problems and for problems with softening behaviors, a subject of research in themselves and no mature tool has emerged from pioneer investigations which can be found, for instance,  in \cite{ladevezemoes1999a, huertarodriguezferran2002, larssonrunesson2003, allixkerfriden2010, PARRETFREAUD.2010.1.1}. Of course, the price that we pay for the expected robustness associated to this choice is a large number of arguably unnecessary degrees of freedom.

\section{Relocalization for the propagation of delamination}\label{sec:reloc}

%\begin{figure}[!t]\centering
%\includegraphics[width=.6\textwidth]{figures/dcb-superstructure.pdf}\caption{Relocalization box around the crack tip}\label{fig:dcb_box}
%\end{figure}

% PK modif
\begin{figure}[!t]\centering
\includegraphics[width=.99\textwidth]{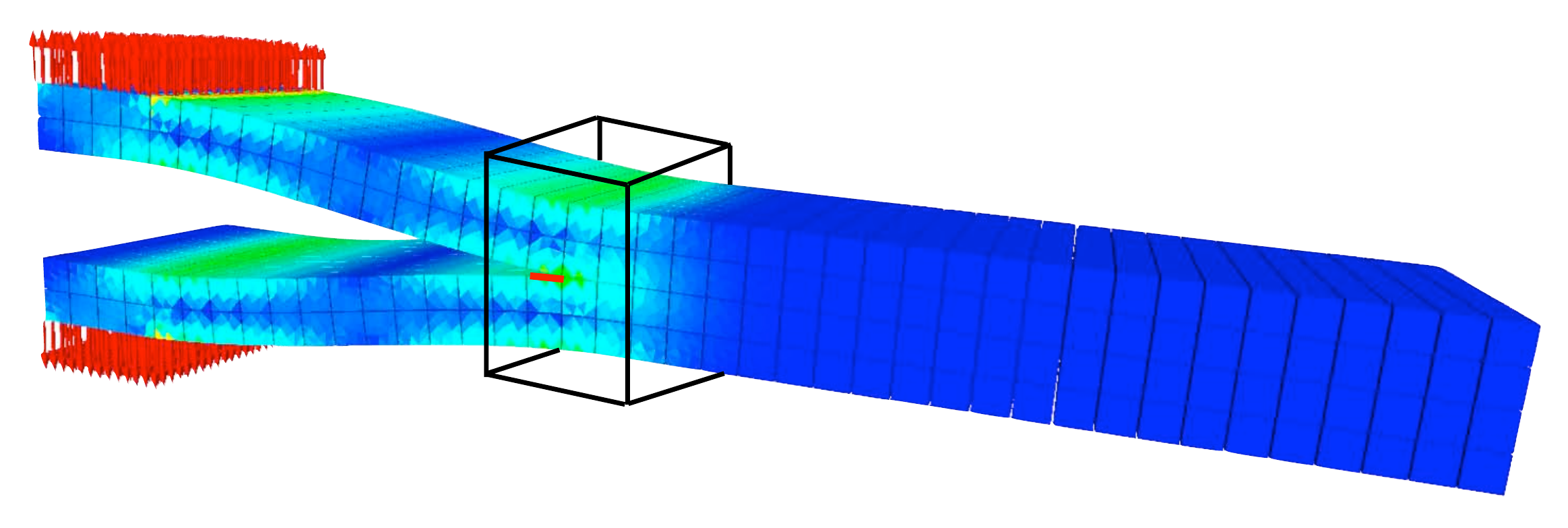}\caption{Relocalization box around the crack tip}\label{fig:dcb_box}
\end{figure}

\begin{figure}[!t]\centering
\begin{minipage}{.18\textwidth}\centering
\includegraphics[width=.99\textwidth]{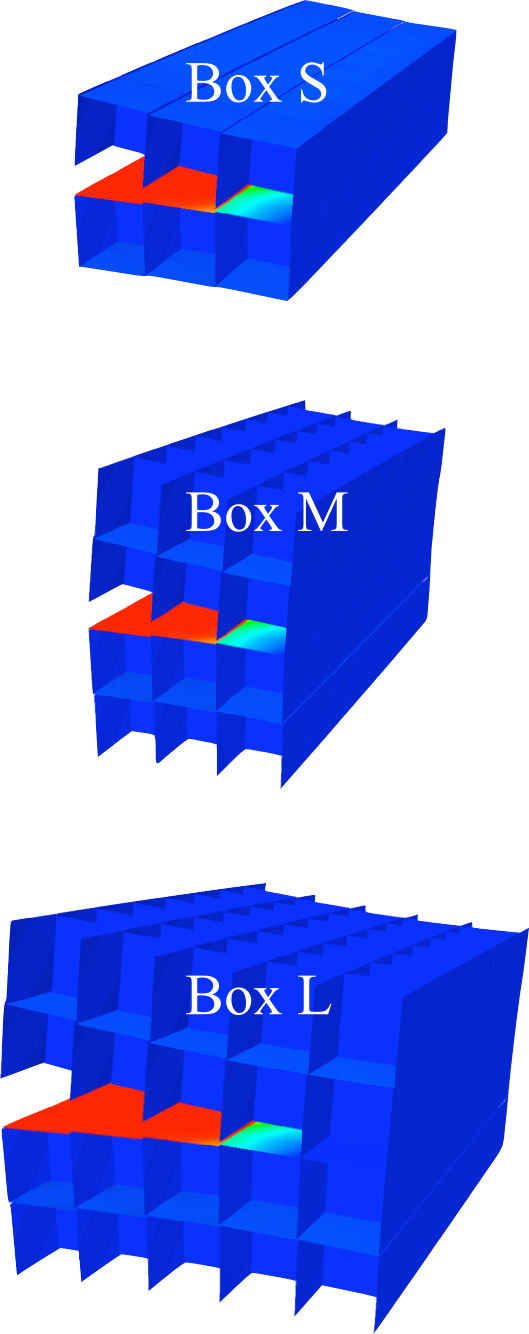}
\end{minipage}
\begin{minipage}{.81\textwidth}\centering
\includegraphics[width=.99\textwidth]{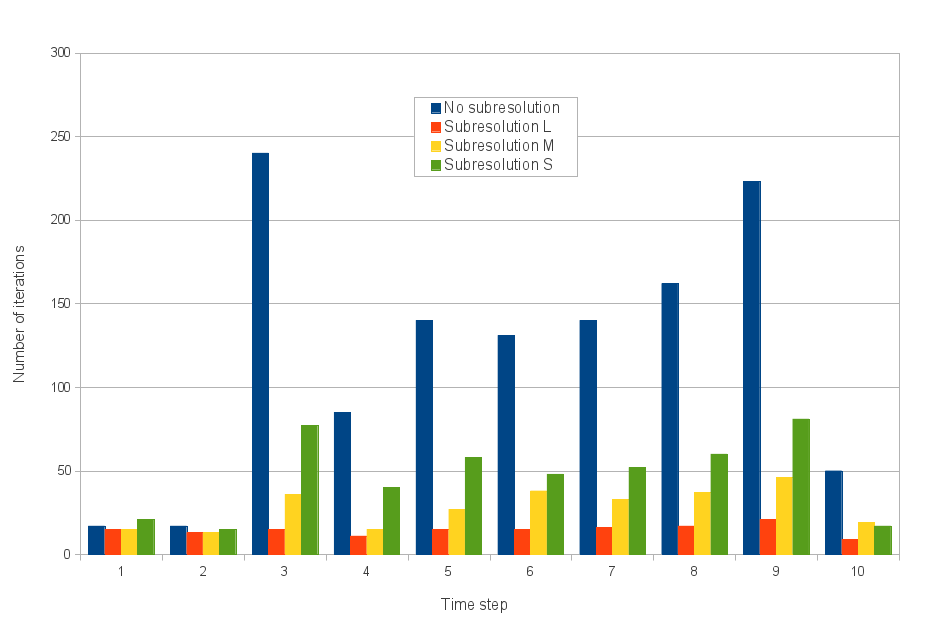}\caption{Performance of the relocalization technique w.r.t. the size of the relocalization box}\label{fig:reloc_box}
\end{minipage}
\end{figure}

Let us consider a simple DCB case, as illustrated in Figure~\ref{fig:dcb_box}. Figure~\ref{fig:reloc_box} shows that the convergence of the classical approach (blue bars) is strongly deteriorated when delamination starts to propagate (time step 3). Indeed, in average 10 times more iterations (and CPU time) are required. This drop in the convergence rate occurring when the cracks propagates can be explained by two driving factors:
\begin{itemize}
\item the singularity near the tip of the crack is very poorly represented by linear macro quantities. Therefore, the complementary part is not localized enough and many iterations are required to transmit them through the structure;
\item the strong non-linearities require many iterations to converge, in particular the prediction of the location of the crack tip at a given time step requires the quasi-convergence of a large number of consecutive equilibrium states.
\end{itemize}
Using non-linear relocalization in a box containing the crack tip and the process zone is a way to deal with those two problems: after each linear stage, a  subproblem is extracted around the crack tip (within the so-called relocalization box whose size is a parameter) as shown in Figure~\ref{fig:dcb_box}. This smaller problem is solved non-linearly with Robin boundary conditions satisfying the macro-equilibrium of the structure.

%Because the non-linearities are strongly localized around the crack tip, non-linear relocalization enables us to deal simultaneously with the two problems. 

Being far enough from the singularity, the macro quantities at the boundary of the box represent the global fields correctly and are then relevant for the prediction of the non-linear evolution of the crack. The resulting modifications only have a small-scale effect which can be transmitted by the upcoming local stage.

The other bars in Figure~\ref{fig:reloc_box} show the influence of the size of the relocalization box on the convergence. Of course the larger the box, the best the filtering of the singularities on the boundary of the box, and the fastest the convergence. The counterpart is that the larger the box, the more expensive the non-linear computations associated to the relocalized solutions, which are obtained by a LaTIn method restricted to the box, are.

\section{Further improvements for delamination/buckling and contact}\label{sec:adaptbuck}
\subsection{Taking into account the structure's slenderness}
Optimal search direction  of one subdomain (for the global step) is known to be the Schur complement of the rest of the structure. In general a very coarse approximation, in practice a scalar, is sufficient to ensure good convergence of the LaTIn method. In \cite{Saavedra12b} it was shown that for plates, the slenderness induces a structural anisotropy that has to be taken into account by the search directions. A dimensional analysis leads to the following relationship between the longitudinal (n) and transverse $(t)$ search directions ($L_\textit0$ is the characteristic length of the structure, $h_\textit0$ is its thickness):
\begin{equation*}
 \frac{\knEo}{\ktEo} =  \left(\frac{L_\textit0}{h_\textit0}\right)^2
\end{equation*}
Here superscript $m$ stands for the micro part of the search direction. The macro part of the search direction has to be set to a large value on perfect interfaces so that the macro continuity of displacements is almost satisfied throughout the iterative process.

\begin{figure}[ht]
       \centering
       \includegraphics[width=.6 \textwidth]{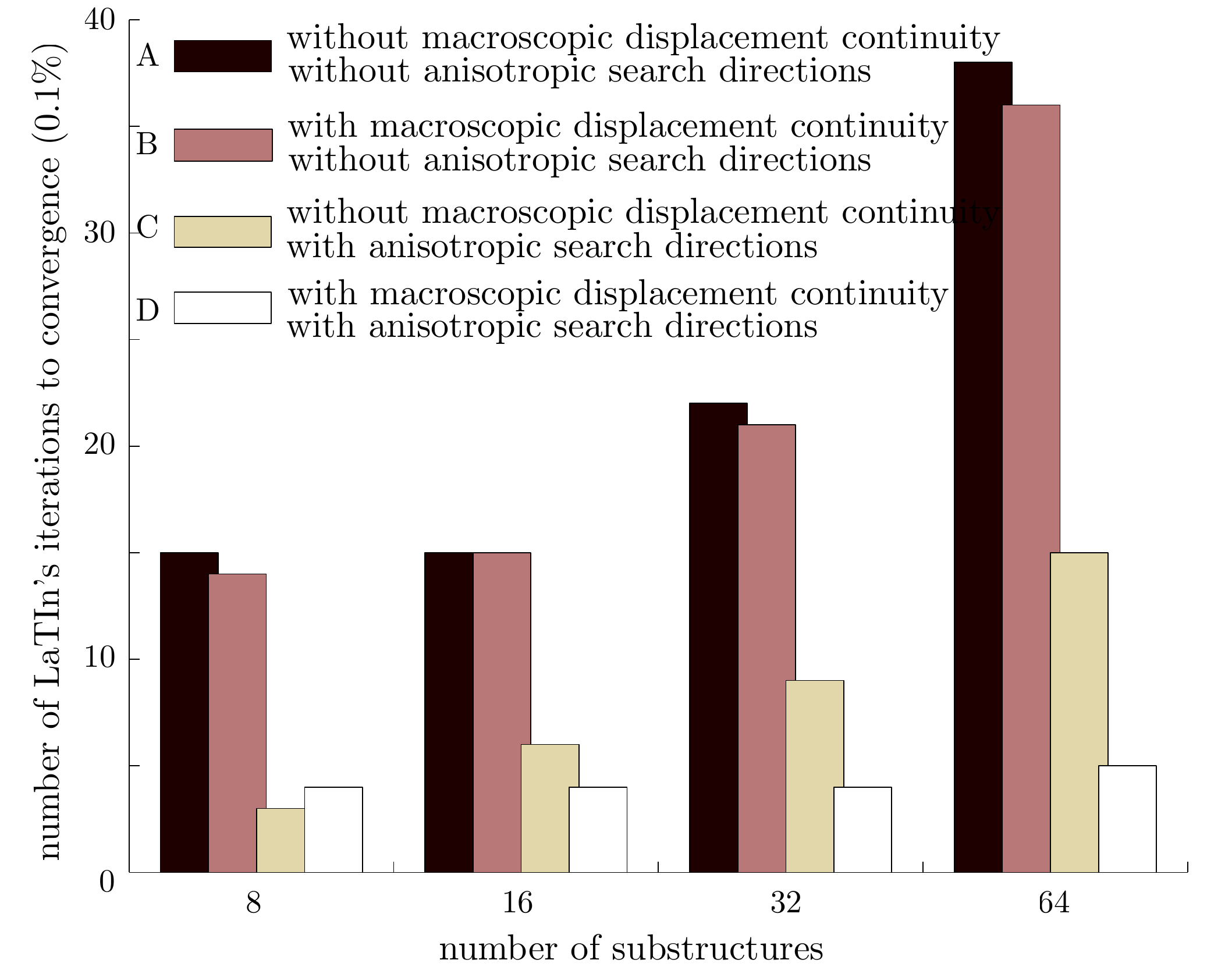}
       \caption{Influence of the search directions and of the substructuration in a bending test case}
       \label{fig:ddr_flexion}
\end{figure}

Figure~\ref{fig:ddr_flexion} shows that these adaptations of the search directions (i.e. anisotropic search directions) are necessary to ensure the scalability of the method.

\subsection{Update of the search directions}
In the context of buckling and contact, the state of the structure can have large influence on its response. This implies that search directions which are optimal in one configuration can be unadapted to others, leading to a drop in the convergence rate. Strategies to update the search directions according to the status of the interfaces must then be settled. Typically for cohesive interfaces, the optimal value is related to the current stiffness of the interface (which depends on the damage state); for fully delaminated interfaces, search directions should be large when contact occurs whereas they should be null when the crack is open.

\begin{figure}[!t]
      \centering
       \includegraphics[width=.8\textwidth]{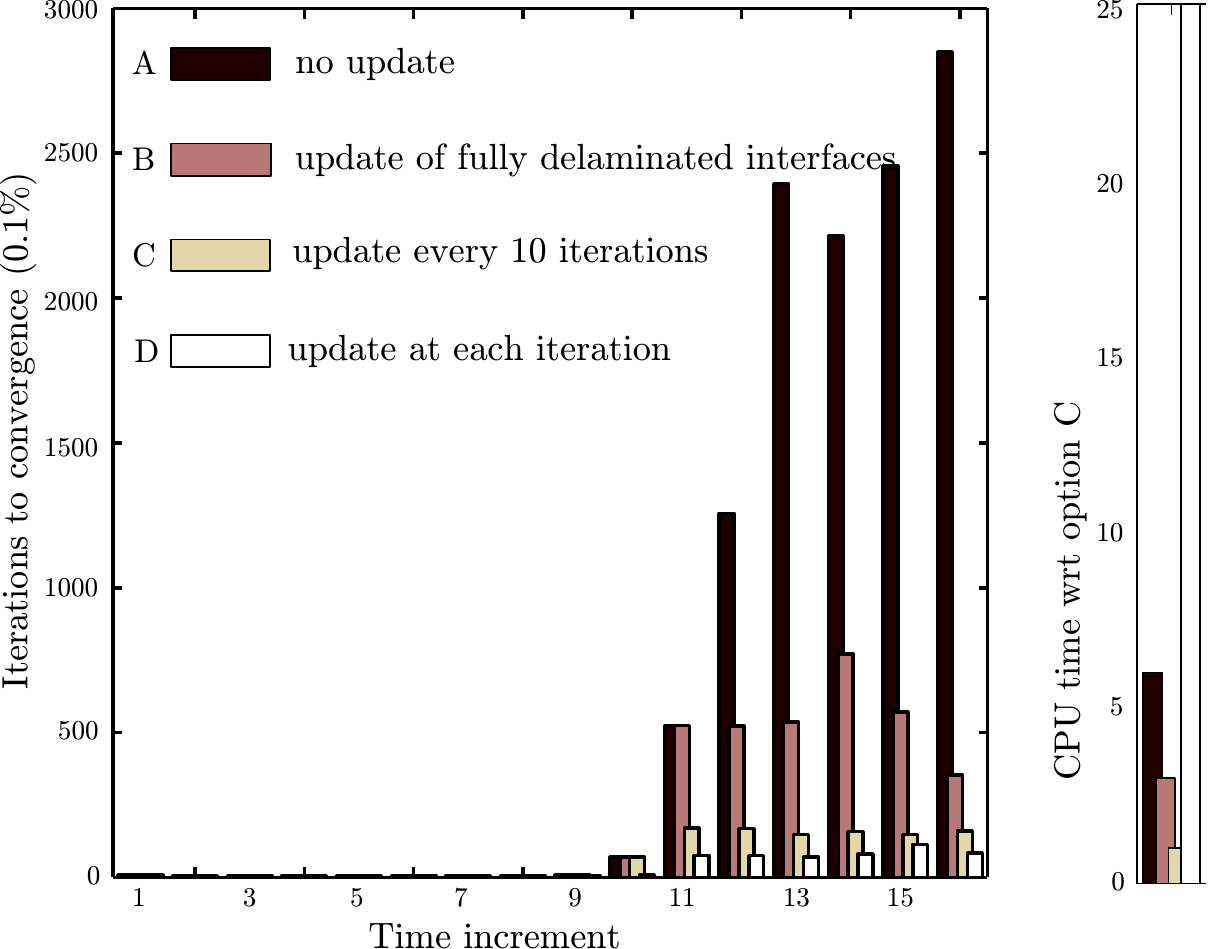}
       \caption{Iterations and CPU time to convergence when using different updating strategies}
       \label{fig:ddr_delam}
\end{figure}

Figure~\ref{fig:ddr_delam} presents the effect of the updating of the search direction according to the damage state. In this simple DCB case, all strategies converge. Of course, the more frequent the updates the less iterations, but each update requires recomputation and factorization of the operators so that a compromise must be found to achieve minimal CPU time.

\begin{figure}[!t]
       \centering
       \includegraphics[width=.7\textwidth]{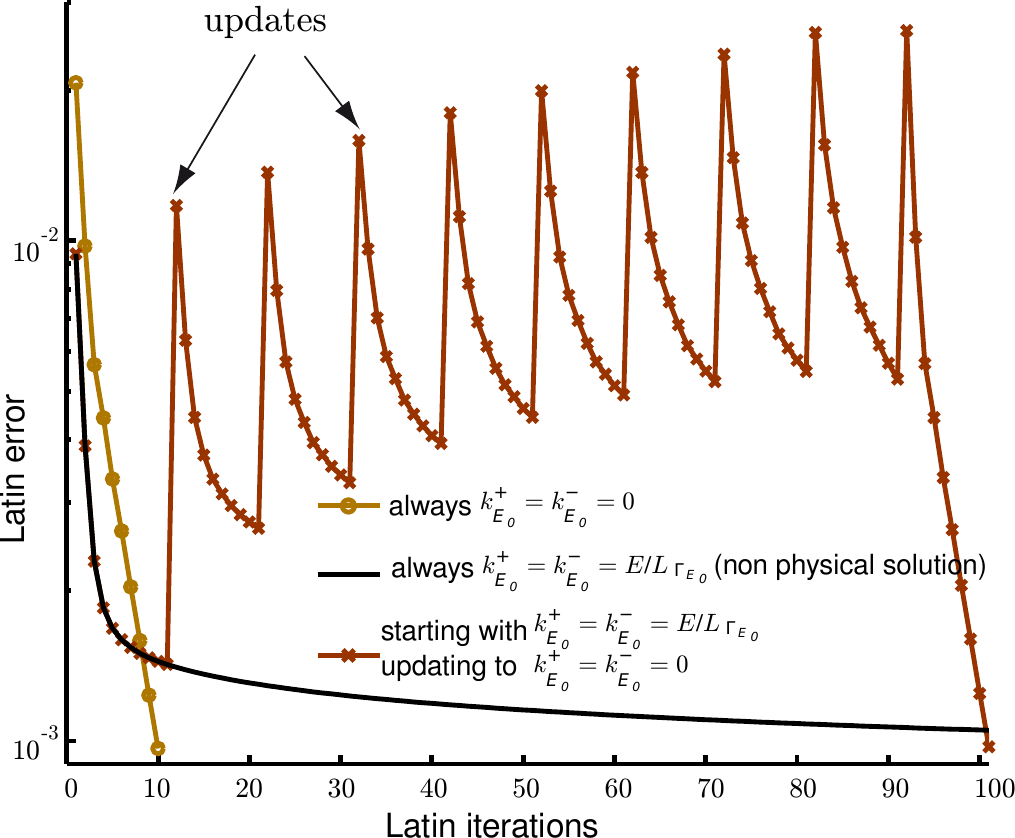}
       \caption{Evolution of the error when using different updating strategies}   \label{fig:contact_open_ddr}
\end{figure}

Figure~\ref{fig:contact_open_ddr} shows the evolution of the error for one time step of a opening contact problem: the plate is initially delaminated in its middle (large closed crack) and a compression loading is imposed. When the optimal search direction is chosen a priori, only 10 iterations are required. When the wrong one is chosen, the algorithm converges to a non-physical solution. Starting from an incorrect initial guess, a simple update procedure (every 10 iterations the parameters are adapted to current contact state) enables to converge to the solution. Note that the convergence is slow because of the unusual extension of the initial contact zone.

\subsection{Pre-delaminated plate sensitive to imperfection} \label{sec:plaque_prefissuree1}

\begin{figure}[t]
       \centering
       \includegraphics[width= \linewidth]{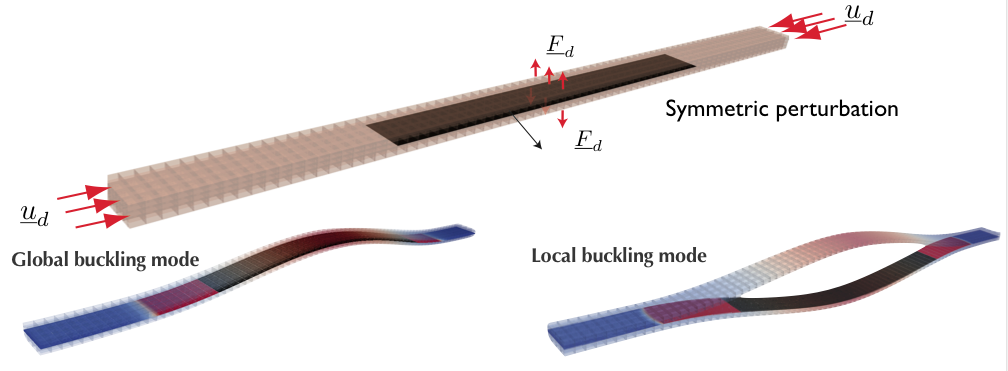}
       \caption{Pre-delaminated plate with associated potential buckling modes}
       \label{fig:pdl-des}
\end{figure}

Let us consider the example of a pre-delaminated plate submitted to compression. The example is designed to be sensitive to imperfection. Indeed two buckling modes can be excited (Figure~\ref{fig:pdl-des}): a global buckling mode (closed interface - mode II delamination) and a local buckling (open interface - mode I delamination). The problem is split into 1280 substructures leading to 3248 interfaces with  about 1.2 millions of d.o.f. The macro-problem involves 30,000 d.o.f. when the super-macro one only involves 132 d.o.f. 

\begin{figure}[!ht]
       \centering
       \includegraphics[width= .6\textwidth]{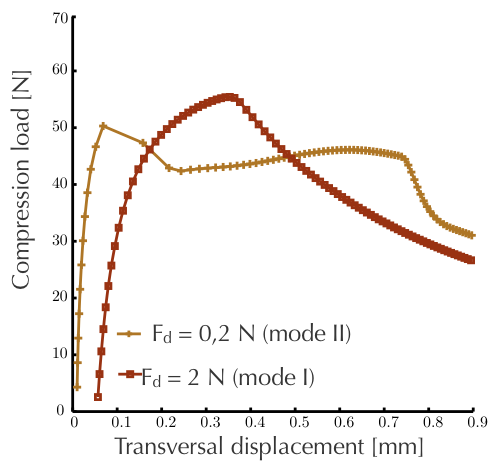}
       \caption{Load displacement curves depending on the imperfection}
       \label{fig:f-w}
\end{figure}

\begin{figure}[!ht]
       \centering
       \includegraphics[width= .85\textwidth]{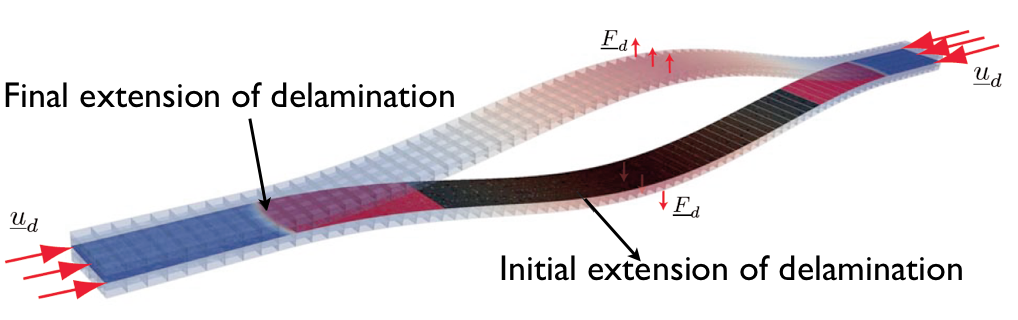}
       \caption{Final state of the plate in case of local buckling}
       \label{fig:final-state}
\end{figure}

%When the imperfection is low the response of the buckling mode of the plate is a global one (mode II delamination, Figure \ref{fig:final-state}), for a displacement of about $0.75$~mm the delamination starts to propagate (Figure~\ref{fig:f-w}- mode II curve).
When the imperfection is small, the response of the buckling mode of the plate is a global one (mode II delamination) and the delamination starts to propagate for a transversal displacement of about $0.75$~mm (see Figure~\ref{fig:f-w} - mode II curve).

When the level of the perturbation is at least ten times greater ($\underline{F}_d=2$~N), a local buckling mode is excited, which leads to the propagation of the delamination for a transversal displacement of about $0.36$~mm (see Figure~\ref{fig:f-w} - mode I curve). Figure \ref{fig:final-state} shows the configuration and the crack propagation of the last time step.

\section{Example of  multiple delamination of a laminated plate in compression}
\label{section:multi_delaminations}
\begin{figure}[ht]
       \centering
       \includegraphics[width= .8\textwidth]{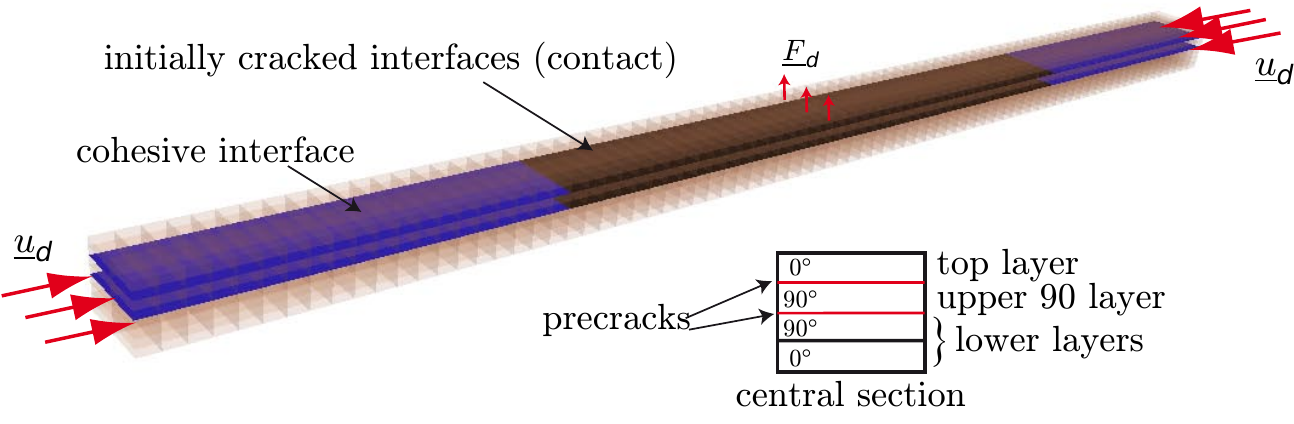}
       \caption{Configuration of the specimen}
       \label{fig:4plis}
\end{figure}
\begin{figure}[!t]
       \centering
       \includegraphics[width= .6\textwidth]{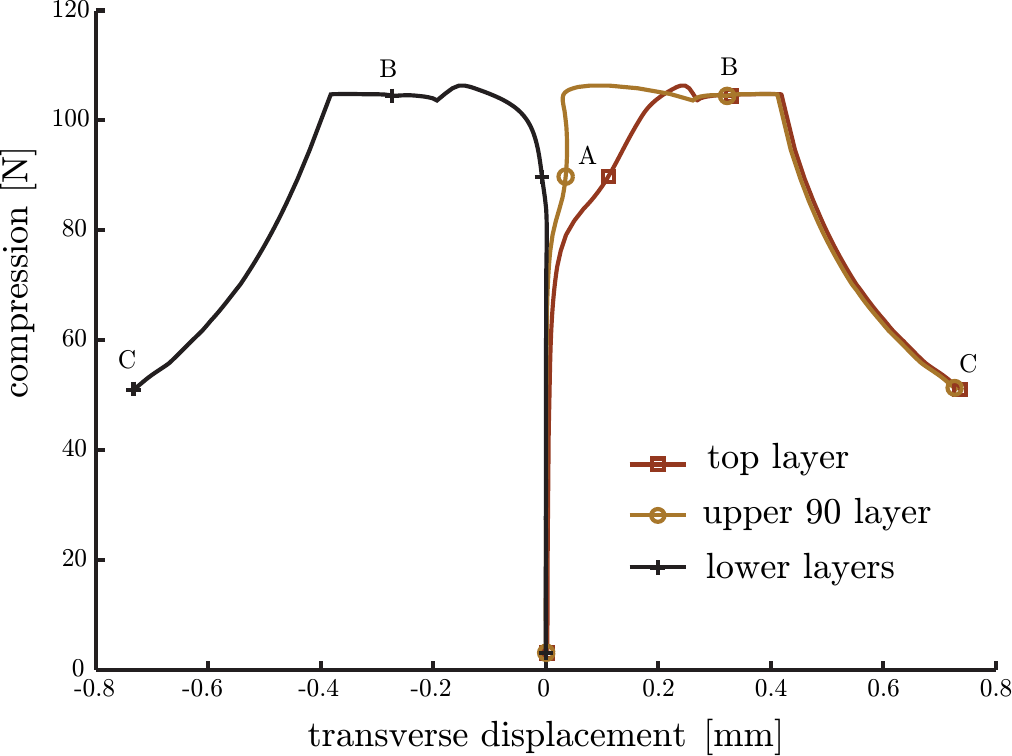}
       \caption{Displacement of three sets of plies during the experiment}
       \label{fig:4plis-F-d}
\end{figure}
\begin{figure}[ht]
       \centering
       \includegraphics[width= .8\textwidth]{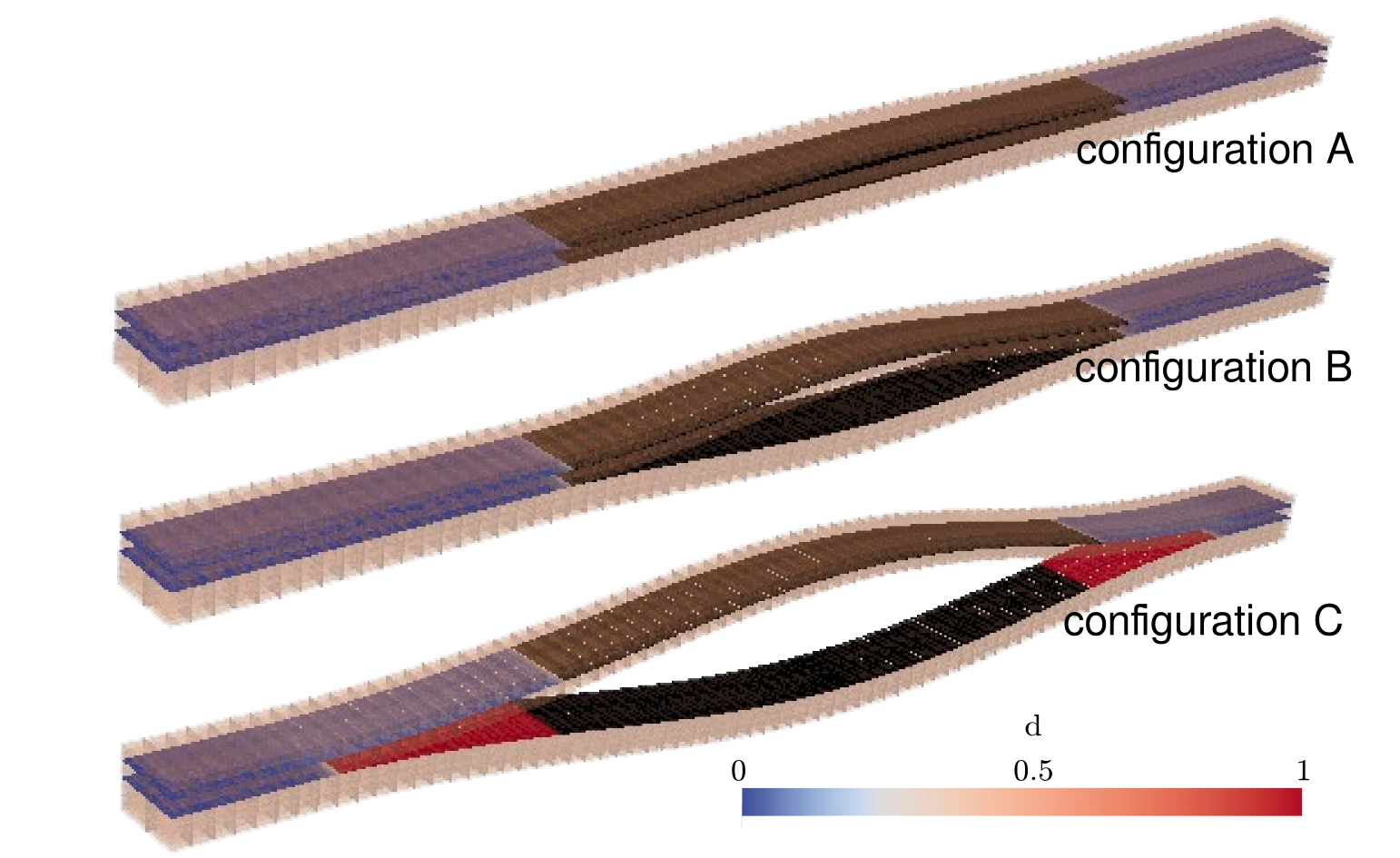}
       \caption{Three configurations. A: all plies separated, B: recontact, C: ruin}
       \label{fig:4plis-def}
\end{figure}

The following example has been chosen because it involves all the difficulties discussed in this last section. Here, the delamination between plies and the possibility of contact and separation between delaminated surfaces in large displacements are studied (others possibles scenarios like transverse cracking are obviously ignored by the model). It is a $[0^\circ/90^\circ]_s$ composite plate, submitted to compression (prescribed displacement $\underline{u}_d$) and  pre-delaminated (with the same initial length) on the central interface and on the upper one. A normal perturbation is applied $\underline{F}_d$ at the level of the top ply (see Figure \ref{fig:4plis}). The layers are assumed to be elastic and orthotropic with the following characteristics in the referenced basis (1- fiber direction, 2- transverse direction in the ply, 3- normal direction):  $E_1=185,500$~MPa, $E_2=E_3=9,900$~MPa, $\nu_{12}=\nu_{13}=0.34$, $\nu_{23}=0.5$, $G_{12}=G_{13}=6,160$~MPa, $G_{23}=3,080$~MPa.  The problem is split into 1,280 substructures leading to 3,248 interfaces, 2 millions of d.o.f. are involved. The macro-problem involves 30,000 d.o.f. when the super-macro one only involves 168 d.o.f. The time interval is discretized into 120 times steps and 30 processors have been used.  
On Figure~\ref{fig:4plis-F-d}, the force-normal displacement responses for the three groups of layers in the middle of the plate are plotted and Figure~\ref{fig:4plis-def} shows three typical configurations. The top ply buckles first  for an applied load of $80$ N leading to a non symmetric configuration (configuration A). In this configuration a loss of contact between the two upper plies and the lower one is observed.  When the applied load reaches $100$ N,  the central and lower plies buckle which leads to a recontact between the upper plies (configuration B). Finally, the configuration C is symmetric. For a value of the transverse displacement of about $0.5$~mm, the delamination starts to propagate in the central interface leading to a decrease of the reaction forces. 



\section{Conclusion}
In order to circumvent the numerous problems associated with the large scale simulation of delamination problems, we have developed a multi-scale framework, wherein a Domain Decomposition approach is optimized to deal with the severe non-linearities. A lot of work is still required to make the approach more generic. The problem of code development is also an issue that needs to be taken into account when developing such dedicated methods. The challenge we would like to address is the one of the Virtual Delamination Testing, in a fully integrated sense. Nevertheless, those studies have been conducted by restricting the potential damage modes of the laminate to delamination. However, it is known that in many cases, like in the case of low-velocity impacts or in the case of the analysis of the tolerance to defects, lumping the damages onto the interfaces cannot reproduce the delamination patterns adequately. In fact, these patterns are often driven by the interaction between transverse cracking and delamination. 

Two main approaches are followed in this more general case. In \cite{Lubineau07} a micro analysis leading to a meso-scale model where a non-local coupling between the interface and the adjacent plies associated with  transverse-cracking / delamination interaction has been developed. A validation of this approach by making use of open-holed tests, which exhibits very different delamination patterns due to scale effect \cite{Green07} can be found in \cite{Abisset11}. The non-locality of the ply model requires to have all the information throughout the thickness of a given ply within a unique sub-domain. This is not a problem with the physical decomposition at the meso-scale proposed in the paper. However, it becomes an issue if one makes use of an automatic decomposition of the structure, as was addressed in \cite{Bordeu10a}. Another approach aims at introducing the transverse cracks explicitly, or at least the transverse cracks that have a strong influence on the delamination pattern \cite{Vandermeer12} by means of X-FEM like technique. In the framework proposed in this paper, following this direction requires to development an hybrid multi-scale  X-FEM method, as was proposed in \cite{Guidault08}.
%such an approach would leads to introduce the multi-scale- Xfem approach \cite{Guidault08} in the framework describe in this paper. 
However, the complexity and highly-intrusive character of such an approach is an issue in terms of code development. 
%Another possibility is to couple, where needed with a non-intrusive method \cite{Gendre09}, the tools presented in this paper with the one developed to simulate a micro-model where all the cracks are introduced explicitly \cite{Violeau09},  as described in \cite{Daghia12}. 
A promising alternative is to perform a local enrichment by means of a non-intrusive coupling of the structural solver with specialized external solvers dedicated to specific failure modes, as was proposed in \cite{Gendre09}. For instance, one can couple the tools presented in this paper with the one developed to simulate a micro-model where all the cracks are introduced explicitly \cite{Violeau09}. This idea was successfully applied in \cite{Daghia12}. 

\bigskip
\noindent\textbf{Acknoledgement}: 
Part of the research leading to these results has received funding from the European Community's Seventh Framework Program FP7/2007-2013 under grant agreement n$^\circ$213371.

\bibliographystyle{plain}
\bibliography{composite-multiechelle}

%\appendix
%\section{Example of appendix}
%Some text here.

\end{document}